\documentclass[12pt]{article}
\newcommand{\ds}{\displaystyle}

\usepackage{graphicx}
\usepackage{latexsym}

\begin{document}

\begin{center}{\bf On Fluid mechanics formulation of 
Monge-Kantorovich Mass Transfer Problem}
\vspace{5mm}

{Kazufumi Ito}
\vspace{3mm}

\noindent Center for Research in Scientific Computation

\noindent North Carolina State University

\noindent Raleigh, North Carolina 27695-8205
\end{center}
\vspace{3mm}

\noindent{\bf Abstract} The Monge-Kantorovich mass transfer problem
is equivalently formulated as an optimal control prblem
for the mass transport equation. The equivalency of
the two problems is establish using the Lax-Hopf formula and 
the optimal control theory arguments.  Also, it is shown that
the optimal solution to the equivalent control problem
is given in a gradient form in terms of the potential
solution to the Monge-Kantorovich problem. It turns out
that the control formulation is a dual formulation 
of the Kantrovich distance problem via the Hamilton-Jacobi equations.
       
\vspace{5mm}




\section{Introduction}

Monge mass transfer problem is that
given two probability density functions $\rho_0(x) \ge 0$ and $\rho_1(x)\ge0$
of $x\in R^d$, find a coordinate map $M$ such that
$$
\int_A \rho_1(x)\,dx=\int_{M(x)\in A}\rho_0(x)\,dx
\leqno (1.1) $$
for all bounded subset $A$ in $R^n$. If $M$ is a smooth one-to-one map,
then it is equivalent to
$$
\mbox{det}(\nabla M)(x)\rho_1(M(x))=\rho_0(x)
\leqno (1.2) $$
where det denotes the determinant of Jacobian matrix of
the map $M$. Clearly, this problem is underdetermined and it
is natural to formulate a costfuctional for the optimal mass transfer.
The so-called Kantorovich (or Wasserstain)
distance between $\rho_0$ and $\rho_1$ is defined by
$$
d(\rho_0,\rho_1)=\inf \int_{R^d} c(x-M(x))\rho_0(x)\,dx.
\leqno (1.3) $$
where $c$ is a convex function and $c(x-y)=c(|x-y|)$ with $c(0)=0$. 
For example $c(x-y)=\frac{1}{p}|x-y|^p$ is 
for the $L^p$ Monge-Kantorovich problem (MKP).
Whenever the infimum is attained by some map $M$, we say that $M$
is an optimal transfer for the Monge-Kantorovich problem.
The Kantorovich distance is the least action that
is necessary to transfer $\rho_0$ into $\rho_1$.

The mass transport problems have attracted a lot of attentions
in recent years and have found applications in many fields of mathematics
such as statistics and fluid mechanics (e.g., see \cite{B2,E}
and \cite{RR} for extensive references). From a more scientific
point of view the Kantorovich distance
provides a valuable quantitative informations to compare
two different density functions and it has been used in various fields
of applications \cite{BB}.

It is shown e.g., in \cite{B1,C,GM,E} that the optimal map $\bar{M}$ is
given by
$$
Dc(x-\bar{M}(x))=\nabla \bar{u}(x)
\leqno (1.4) $$
for a potential function $\bar{u}$, where $Dc$ denotes the derivative of $c$. 
In fact $\bar{u}$ is the optimal solution to the Kantorovich dual 
problem (2.2).
If $c$ is uniformly convex, then
we can solve (1.4) for $\bar{M}$ in terms of $\nabla \bar{u}$.  
For example for $L^p$ MKP
$$
x-\bar{M}(x)=|\nabla\bar{u}(x)|^{q-2}\,\nabla\bar{u}(x) \mbox{  with  }
\frac{1}{p}+\frac{1}{q}=1.
$$
For $L^2$  MKP, it follows from (1.2) and (1.4) that if 
$\psi=\frac{|x|^2}{2}-\bar{u}(x)$, 
then $\bar{M}(x)=x-\nabla \bar{u}=\nabla\psi$
and thus $\psi$ satisfies the Monge-Ampere equation
$$
\mbox{det}(H \psi)(x)\rho_1(\nabla \psi )=\rho_0(x),
\leqno (1.5) $$
where $H\psi$ is the Hessian of $\psi$. 

In \cite{BB} the $L^2$ MKP is equivalently reformulated as an optimal 
control problem:
$$
d(\rho_0,\rho_1)= \min\;\; \left[\frac{1}{2}\int^1_0\int_{R^d} 
\rho(t,x)\, |V(t,x)|^2\,dxdt \quad\mbox{over vector field $V=V(t,x)$}\right]
\leqno (1.6) $$
subject to
$$\begin{array}{l}
\rho_t+\nabla \cdot(\rho\,V)=0
\\ \\
\rho(0,x)=\rho_0(x) \quad \mbox{and}\quad \rho(1,x)=\rho_1(x)
\end{array} \leqno (1.7) $$
Moreover if $\bar{V}(t,x)$ is an optimal solution to (1.6)--(1.7) and 
the Lagrange coordinate $\bar{X}(t;x)$ satisfies
$$
\frac{d}{dt}\bar{X}=\bar{V}(t,\bar{X}(t;x)),\quad \bar{X}(0;x)=x,
$$
then $\bar{M}(x)=\bar{X}(T;x)$.

The contribution of this paper is that 
we will show that the optimal vector field $\bar{V}$
to problem (1.6)--(1.7) is given by
$$
\bar{V}(x,t)=\nabla_x\bar{\phi}(t,x)
\leqno (1.8) $$
where the potential function $\bar{\phi}$ satisfies 
the Hamilton-Jacobi equation
$$
\bar{\phi}_t+\frac{1}{2}\,|\nabla \bar{\phi}|^2=0,
\quad \bar{\phi}(0,x)=-\bar{u}(x)
\leqno (1.9) $$
and  $\bar{u}$ determines the optimal map $\bar{M}$ in (1.4).  
Thus, (1.8)--(1.9) is an optimal feedback solution to 
control problem (1.6)--(1.7), i.e., given $\rho_0,\;\rho_1$
first we determine $\psi$ by (1.5) and let $\ds \bar{u}=\frac{|x|^2}{2}-\psi$
and then determine $\bar{V}$ by (1.8)--(1.9). 
\newpage

\noindent Moreover, it will be shown that
$$
d(\rho_0,\rho_1)=\min\;\;\left[\int_{R^d} (\rho_1(x)v(x)-\rho_0(x)
\phi(0,x))\,dx \mbox{  over  } v \right]
$$
subject to 
$$
\phi_t+\frac{1}{2}\,|\nabla \phi|^2=0,\quad \phi(1,x)=v(x).
$$
It is the other control formulation of the $L^2$ MKP and
is an optimization problem over the potential fuction $v$
subject to the Hamilton-Jacobi equation.  

For the non-quadratic $c$ case, 
the (generalized) optimal control problem is formulated as
$$
\min\quad \int^1_0\int_{R^d} \rho(t,x)\,c(V(t,x))\,dxdt
\leqno (1.10) $$
subject to (1.7). 
In this case the optimal vector field $\bar{V}$ is given by
$$
\bar{V}=Dc^*(\nabla \bar{\phi})
\leqno (1.11) $$
where $c^*$ is the convex conjugate function of $c$ defined by
$$
c^*(y)=\sup_x\; \{x\cdot y-c(x)\}.
$$
For $L^p$ MKP
$$
c(x)=\frac{1}{p}|x|^p,\;x\in R^d,\quad c^*(y)=\frac{1}{q}|y|^q,\;\;y\in R^d,  
$$
and
$$
\bar{V}(t,x)=|\nabla_x\phi(t,x)|^{q-2}\,\nabla_x\phi(t,x)
$$
where $\ds \frac{1}{p}+\frac{1}{q}=1$ and $p\in (1,\infty)$.  
The potential function $\phi=\phi(t,x)$ satisfies
$$
\bar{\phi}_t+c^*(\nabla \bar{\phi})=0,\quad \bar{\phi}(0,x)=-\bar{u}(x).
\leqno (1.12) $$
If $c$ is uniformly convex, then from (1.4) 
$$
\bar{M}(x)=x-Dc^*(\nabla \bar{u}(x)).
$$
Thus, from (1.2) $\bar{u}$ satisfies
$$
\mbox{det}(\nabla \bar{M})(x)\rho_1(x-Dc^*(\nabla \bar{u}(x)))=\rho_0(x).
\leqno (1.13) $$
For $L^2$ MKP (1.13) is reduced to (1.5).
Hence the optimal solution to (1.10) subject to (1.7) is given
in the feedback form (1.11)-(1.13). 

An outline of our presentation is as follows.
In Section 2 the basic theoretical results concerning
the MKP problem is reviewed following \cite{E}. Then equivalent
variational formulations (2.5) and (2.8) for the potential function  
are then derived using the duality and the
Lax-Hoph formula. In Section 3 we present formal arguments
that show the feedback solution (1.7)--(1.9) to (1.5)--(1.6).
In Section 4 we validate the steps in Section 3 mathmatically
for $L^2$ MKP. In Section 5 we present the proofs for the general case.  

\section{Variational Formulations}

In order to present our treatment of the MKP problem,
we first recall a basic theoretical result in this section.
The following relaxed problem of (1.3) is introduced by Kantorovich.
Let $\cal{M}$ be a class of random probability measures $\mu$
on $R^d\times R^d$ satisfying proj$_y\mu=\rho_0\,dx$
and proj$_x\mu=\rho_1\,dy$. Then we define the relaxed cost-functional
$$
J(\mu)=\int_{R^d\times R^d} c(x-y)\,d\mu(x,y) \quad
\mbox{over}\;\; \cal{M}.
\leqno (2.1) $$
Consider the dual problem of (2.1); maximize
$$\begin{array}{l}
\ds \int_{R^d} u(x)\rho_0(x)\,dx+\int_{R^d} v(x)\rho_1(y)\,dy
\\ \\
\ds \quad \mbox{subject to}\;\;
u(x)+v(y) \le c(x-y).
\end{array} \leqno (2.2) $$
The point of course is that the Lagrange multiplier associated with
the inequality in (2.2) solves problem (2.1).
The following theorem \cite{B1,E,C,GM} provides the solution to 
(2.2) and (1.3).
\vspace{10mm}

\noindent{\bf Theorem 2.1} 

\noindent (1) there exists a maximizer $(\bar{u},\bar{v})$ of
problem (2.2).
\vspace{2mm}

\noindent (2) $(\bar{u},\bar{v})$ are dual $c$-conjugate functions, i.e.,
$$\begin{array}{l}
\ds \bar{u}(x)=\inf_y\;(c(x-y)-\bar{v}(y))
\\ \\
\ds \bar{v}(y)=\inf_x\;(c(x-y)-\bar{u}(x))
\end{array} $$

\noindent (3) $\bar{M}(x)$ satisfying 
$Dc(x-\bar{M}(x))=\nabla\bar{u}(x)$  solves MKP problem.
\vspace{2mm}

It follows from Theorem 2.1 that (2.3) is reduced to maximizing
$$
J(u) = \int_{R^d} u(x)\rho_0(x)\,dx+\int_{R^d} v(y)\rho_1(y)\,dy
\leqno (2.4) $$
over functions $u$, where $v$ is the $c$-conjugate function 
of $u$. The $c$-conjugate function of a function $u$ is defined by 
$$
v(y)=\inf_x\;(c(x-y)-u(x)).
$$
It is easy to show that the bi $c$-conjugate function $\tilde{u}$ of $u$ 
satisfies $\tilde{u}\ge u$ a.e. and thus the maximizing pair 
$(u,v)$ of (2.4) is automatically $c$-conjugate each other.
Similarly, we have the equivalent problem of maximizing
$$
J(v)=\int_{R^d} u(x)\rho_0(x)\,dx+\int_{R^d} v(y)\rho_1(y)\,dy
\leqno (2.5) $$
where
$$
u(x)=\inf_x\;(c(x-y)-v(y)).
$$

Let $c^*$ be the convex conjugate of $c$, i.e.,
$$
c^*(x)=\sup_y((x,y)-c(y)). 
$$
By the Lax-Hopf formula \cite{E1}, if $\phi$ is the viscosity solution to
$$
\phi_t+c^*(\nabla\phi)=0,\quad \phi(1,y)=v(y)
\leqno (2.6) $$
then 
$$
\phi(0,x)=\sup_y\;(v(y)-c(x-y))=-u(x).
\leqno (2.7) $$
Thus, Problem (2.2) can be equivalently formulated as maximizing
$$
J(v)=\int_{R^n} (\rho_1(x)v(x)-\rho_0(x)\phi(0,x)) \,dx
\leqno (2.8) $$
subject to (2.6). 

\section{Derivation of Optimal Feedback Solution}

The optimality condition of (2.8) subject (2.6) 
is formally derived as follows. We define the Lagrangian
$$
L(\phi,\lambda)=J(\phi(1))-\int^1_0 \int_{R^d} 
(\phi_t+c^*(\nabla\phi))\lambda\,dxdt.
\leqno (3.1) $$
By applying the Lagrange multiplier theory
the necessary optimality  is given by
$$\begin{array}{l}
\ds L_\phi(\phi,\lambda)(h)=\int^1_0\int_{R^d} (\lambda_t
+(Dc^*(\nabla\phi)\,\lambda)_x)
h\,dxdt
\\ \\
\ds\;\;  -\int_{R^d}(h(1,x)\lambda(1,x)-h(0,x)\lambda(0,x))\,dx
+\int_{R^d}(h(1,x)\rho_1(x)-h(0,x)\rho_0(x))\,dx=0
\end{array} \leqno (3.2) $$
for all $h\in C^1_0([0,1] \times R^d)$. 
Hence the necessary optimality reduces to
$$\begin{array}{l}
\lambda_t+(Dc^*(\nabla\bar{\phi})\,\lambda)_x=0
\\ \\
\lambda(0,x)=\rho_0(x),\quad \lambda(1,x)=\rho_1(x).
\end{array} \leqno (3.3) $$
This implies that if we let $\bar{V}(t,x)=Dc^*(\nabla\bar{\phi}(t,x))$ in 
$$
\bar{\rho}_t+(\bar{V}\,\bar{\rho})_x=0,\quad \bar{\rho}(0,x)=\rho_0,
$$
then $\bar{\rho}(1,x)=\rho_1(x)$. Moreover, we can argue that
$$
\int_{R^d} (\rho_0 \bar{u}(x)+\rho_1(x)\bar{v}(x))\,dx
-\int^1_0 \int_{R^d} \bar{\rho}(t,x)c(\bar{V}(t,x))\,dxdt=0
\leqno (3.4) $$   
since 
$$
c(\bar{V})=(\nabla_x \bar{\phi},Dc^*(\nabla_x\bar{\phi}))
-c^*(\nabla_x\bar{\phi}).
$$
It follows from (3.3)--(3.4) that 
$\bar{V}=Dc^*(\nabla\bar{\phi})$ is the optimal solution to
(1.10) subject to (1.6). In fact, for sufficiently smooth pair $(\rho,V)$ 
satisfying (1.6), we define Lagrange coordinate $X(t,x)$ by
$$
\frac{d}{dt}X=V(t,X(t,x)),\quad X(0,x)=x,
$$
Then for all test function $f$
$$\begin{array}{l}
\ds \int^1_0\int_{R^d} f(t,x)\rho(t,x)\,dxdt=
\int^1_0\int_{R^d} f(t,X(t,x))\rho_0(x)\,dxdt.
\\ \\
\ds \int^1_0\int_{R^d} f(t,x)\rho(t,x)V(t,x)\,dxdt
=\int^1_0\int_{R^d} V(t,x(t))f(t,X(t,x))\rho_0(x)\,dxdt.
\end{array} \leqno (3.5) $$
Note that (1.6) and (3.5) imply that $M(x)=X(1,x)$ satisfies condition (1.1).
Letting $f=c(V)$ in (3.5), we have
$$\begin{array}{l}
\ds I=\int^1_0\int_{R^d}\rho(t,x)c(V(t,x))\,dxdt= 
\int^1_0\int_{R^d} c(V(t,X(t,x))\rho_0(x)\,dxdt.
\\ \\
\ds\qquad =\int^1_0\int_{R^d} c(\frac{d}{dt}X(t,x))\rho_0(x)\,dxdt \ge
\int c(X(1,x)-X(0,x))\rho_0(x)\,dx.
\end{array} $$
where we used the Jessen's inequality.
Since $\bar{M}(x)$ is the optimal solution to 
(1.3), it follows that
$$
I \ge \int_{R^d} c(x-\bar{M}(x))\rho_0(x)\,dx=d(\rho_0,\rho_1).
\leqno (3.6) $$
From (3.4), (3.6) and Theorem 2.1 
$$
d(\rho_0,\rho_1)=\int_{R^d} (\rho_0 \bar{u}(x)+\rho_1(x)\bar{v}(x))\,dx
=\int^1_0 \int_{R^d} \bar{\rho}(t,x)c(\bar{V}(t,x))\,dxdt \le I.
\leqno (3.7) $$
for all pair $(\rho,V)$ satisfying (1.6). 
That is, $(\bar{\rho},\bar{V})$ is optimal.

\section{Proof of (3.3)--(3.4)}

In this section we give a proof for the steps of deriving
the optimality condition (3.3) and equality (3.4)
in the case when $p=2$, i.e., $c(|x-y|)=\frac{1}{2}\,|x-y|^2$.
Suppose $v\in W^{1,\infty}(R^d)$ and $v$ is semi-convex. 
Then it follows from the Lax-Hopf formula 
$$
\phi(t,x)=\sup_y\;\{v(y)-(1-t)c(\frac{x-y}{1-t})\}
$$
(e.g., see \cite{E1}) that (2.6) has a unique solution 
$\phi \in W^{1,\infty}([0,1] \times R^d)$  with 
$$
|\phi(t)|_{W^{1,\infty}} \le |v|_{W^{1,\infty}} \quad
\mbox{and}\quad |\phi_t(t)|_{L^\infty} \le \frac{1}{2}|v|_{W^{1,\infty}}^2
\leqno (4.1) $$
and 
$$
\phi(t,x+z)-2\phi(t,x)+\phi(t,x-z) \ge
-C\,|z|^2  \mbox{  for all $t\in [0,1]$ and $x,z\in R^d$}. 
\leqno (4.2) $$
where we assumed $\ds v+\frac{C}{2}|x|^2$ is convex. 
Let $\phi^\tau$ be the solution to (2.6) with $\phi^\tau(1)=v+\tau\,h$
for $h \in C^2_0(R^d)$. Assume $\ds \phi^\tau(1)
+\frac{C}{2}|x|^2$ be convex for $|\tau|\le 1$ and thus (4.2)
holds for $\phi^\tau$. 
\vspace{2mm}

\noindent\underline{Step 1} Since
$y \to c(x-y)-v(y)$ is coersive, for each $x \in R^d$ there exist  
$y,\;y^\tau\in R^d$ such that
$$
\phi(0,x)=v(y)-c(x-y),\quad \phi^\tau(0,x)=(v+\tau\,h)(y^\tau)
-c(x-y^\tau).
$$ 
Thus
$$
\phi(0,x)\ge v(y^\tau)-c(x-y^\tau)=(x+\tau\,h)(y^\tau)
-c(x-y^\tau)+\tau\,h(y^\tau)
$$
and
$$
\phi(0,x)-\phi^\tau(0,x)\ge \tau\,h(y^\tau)
$$
Similarly
$$
\phi^\tau(0,x)-\phi(0,x)\ge \tau\,h(y).
$$
Hence
$$
|\phi^\tau(0,\cdot)-\phi(0,\cdot)|_\infty \le \tau\,|h|_\infty.
\leqno (4.3) $$

\noindent\underline{Step 2} Note that
$$
(\phi^\tau-\phi)_t+
\frac{1}{2}(\nabla\phi^\tau+\nabla\phi)\cdot
(\nabla\phi^\tau-\nabla\phi)=0.
\leqno (4.4) $$ 
Let $\eta_\epsilon,\;\epsilon>0$ be the standard molifier.
Then 
$$
|\nabla(\eta_\epsilon*\phi)|_\infty\le |\nabla \phi|_\infty
\leqno (4.5) $$
and
$$
\nabla(\eta_\epsilon*\phi) \to \nabla \phi
\quad\mbox{a.e. as $\epsilon\to 0^+$}.
\leqno (4.6) $$
Moreover (4.2) implies
$$
D^2(\eta_\epsilon *\phi) \ge -C.
$$  
Thus
$$
(\nabla (\eta_\epsilon *\phi),\nabla\psi) \le dC \int_{R^d}\psi\,dx
$$
for $\psi \in W^{1,1}(R^d)$ and $\psi\ge0$ a.e. in $R^d$.   
Thus from (4.5)--(4.6) and the Lebesgue dominated convergence theorem,
letting $\epsilon\to0^+$
$$
(\nabla \phi,\nabla\psi) \le dC \int_{R^d}\psi\,dx
\leqno (4.7) $$
for $\psi \in W^{1,1}(R^d)$ and $\psi\ge0$ a.e. in $R^d$.
It now follows from (4.4) and (4.7) that  
$$
\frac{d}{dt}|\phi^\tau(t,\cdot)-\phi(t,\cdot)|_1
\ge -dC \,|\phi^\tau(t,\cdot)-\phi(t,\cdot)|_1,\quad
\phi^\tau(1)=\phi(1)+\tau\,h
$$
and thus
$$
|\phi^\tau(0,\cdot)-\phi(0,\cdot)|_1 \le \tau\,e^{dC}\,|h|_1.
\leqno (4.8) $$
Since from (2.6) 
$$
\int^1_0 \int_{R^d} \frac{1}{2}|\nabla\phi^\tau|^2 \, dxdt  
=\int_{R^d}(\phi^\tau(0,x)-v(x))\,dx,
$$
we have 
$$
\int^1_0 \int_{R^d} |\nabla\phi^\tau|^2 \, dxdt \to 
\int^1_0 \int_{R^d} |\nabla\phi|^2 \, dxdt.
$$
as $\tau\to0$. Since $L^2((0,1)\times R^d)$ is a Hilbert space, 
this implies that 
$$
\int^1_0 \int_{R^d} |\nabla\phi^\tau-\nabla\phi|^2 \, dxdt \to 0 
\leqno (4.9) $$
as $\tau \to 0$.
\vspace{2mm}

\noindent\underline{Step 3} For $\epsilon >0$ let us consider 
$$
\lambda_t+(\nabla\phi\,\lambda)_x=\epsilon\,\Delta \lambda, 
\quad \lambda(0)=\rho_0
$$
Since $\phi$ is Lipschitz on $[0,1]\times R^d$,
$$
t\to \int_{R^d} [\epsilon\,(\nabla \lambda,\nabla \psi)
-(\nabla\phi(t,\cdot)\lambda,\nabla\psi)]\,dx
$$
defines an integrable, bounded, coersive form on $H^1(R^d)\times H^1(R^d)$
and thus it follows from the parabolic equation theory 
(e.g., see \cite{Ta,IK}) that there exits a unique solution 
$\lambda_\epsilon \in H^1(0,1;L^2(R^d)) \cap L^2(0,1;H^2(R^d))$
provided that $\rho_0 \in H^1(R^d) \cap L^1(R^d) \cap L^\infty(R^d)$. Moreover
$$\begin{array}{l}
\ds |\lambda_\epsilon(t)|_1 \le |\rho_0|_1,
\\ \\
\ds \frac{1}{2}\,(|\lambda_\epsilon(1)|_2^2-|\rho_0|_2^2)
+\epsilon \int^1_0|\nabla \lambda_\epsilon(t)|^2_2\,dt=0
\\ \\ 
|\lambda_\epsilon(t)|_\infty \le e^{Ct}\,|\rho_0|_\infty.
\end{array} \leqno (4.10) $$ 
For the last estimate we have from (4.7)
$$
\frac{1}{p}\,\frac{d}{dt}\,|\lambda_\epsilon|^p \le 
\frac{p-1}{p}(\nabla\phi,\nabla |\lambda_\epsilon|^p) 
\le \frac{(p-1)dC}{p}\,|\lambda_\epsilon|^p.
$$
for $p\ge 1$ and thus
$$
|\lambda_\epsilon|_p \le e^{\frac{p-1}{p}dCt}\,|\rho_0|_p.
$$ 
Thus $\lambda_\epsilon$ is uniformly bounded in 
$L^2((0,1)\times R^d)$.
Hence there exists a $\lambda\in L^\infty (0,1;L^1(R^d) \cap L^\infty(R^d))$
and subsequaence of $\lambda_\epsilon$ (denoted by the same)
such that $\lambda_\epsilon$ converges weakly to $\lambda$ 
in $L^2((0,1)\times R^d)$ and $\lambda_\epsilon(1) \to \lambda(1)$
in $L^2(\Omega)$. Since for $\psi \in C^1_0([0,1] \times R^d)$
$$
\int^1_0\int_{R^d} (\lambda_\epsilon \psi_t+\lambda_\epsilon 
\nabla \phi\cdot\nabla \psi
-\epsilon\,\nabla \lambda_\epsilon\cdot\nabla\psi)\,dxdt
=\int_{R^d} (\rho_0\psi(0)-\lambda_\epsilon(1)\psi(1))\,dx
\leqno (4.11) $$
it follows from (4.10)--(4.11) that letting $\epsilon \to 0^+$
$$
\int^1_0\int_{R^d} (\lambda\psi_t+\lambda \nabla \phi\cdot\nabla \psi)\,dxdt
=\int_{R^d} (\rho_0\psi(0)-\lambda(1)\psi(1))\,dx.
\leqno (4.12) $$
Hence $\lambda$ is a weak solution to 
$$
\lambda_t+(\lambda\,\nabla \phi)=0,\quad \lambda(0)=\rho_0.
\leqno (4.13) $$

Next we show that (4.13) has the weak unique solution
in $L^\infty(0,1;L^1(R^d)\cap L^1(R^d))$.
Let $\eta_\epsilon,\;\epsilon>0$ be the standard molifier
and consider the adjoint equation 
$$
\psi_t+\nabla(\eta_\epsilon*\phi)\cdot\nabla\psi=
f \in C^\infty_0([0,1]\times R^d),\quad\psi(1)=0.
\leqno (4.14) $$
Then, (4.14) has a smooth unique sulution $\psi$ and
$J=|\nabla \psi|$ satisfies
$$
J_t+ \nabla(\eta_\epsilon*\phi)\cdot\nabla J+
D^2(\eta_\epsilon*\phi)J=\nabla f, \quad J(1)=0.
\leqno (4.15) $$
Since $\psi$ has compact support, $J$ has a positive 
maximum over $[0,1]\times R^d$ at some point $(t_0,x_0)$. 
If $0\le t_0 <1$, then from (4.15)
$$
J_t(t_0,x_0) \le 0 \mbox{  and  } \nabla J (t_0,x_0)= 0.
$$
Thus,
$$
D^2(\eta_\epsilon*\phi)J(t_0,x_0) \ge \nabla f (t_0,x_0)
$$
Since from (4.2) $D^2(\eta_\epsilon*\phi)\le -C$, this implies
$$
|\nabla \psi|_\infty=J(t_0,x_0) \le \frac{|\nabla f|_\infty}{C}
\leqno (4.16) $$
Let $\lambda,\;\tilde{\lambda}$ is two weak solutions to (4.13).
Then, it folows from (4.12) and (4.14) that
$$
\int^1_0\int_{R^d}(\lambda-\tilde{\lambda})f\,dxdt
=\int^1_0\int_{R^d}(\lambda-\tilde{\lambda})
(\nabla(\eta_\epsilon *\phi)-\nabla\phi)\nabla\psi\,dxdt.
$$
By letting $\epsilon\to0^+$, it follows from (4.5)--(4.6), (4.16) 
and the Lebesgue dominated convergence theorem that
$$
\int^1_0\int_{R^d}(\lambda-\tilde{\lambda})f\,dxdt=0
$$
for all $f \in C^\infty_0([0,1]\times R^d)$ and therfore
$\lambda=\tilde{\lambda}$.

Now, let $\lambda^\tau$ be the solution to (4.13) associated with $\phi^\tau$.
Since $\lambda^\tau$ is uniformly bounded in $L^\infty (0,1;L^1(R^d)
\cap L^\infty(R^d))$,   
there exists a $\lambda^* \in L^\infty (0,1;L^1(R^d)\cap L^\infty(R^d))$
such that $\lambda^\tau$ converges weakly to $\lambda^*$ in 
$L^2((0,1)\times R^d)$ and $\lambda^\tau (1)$ converges weakly to 
$\lambda^*(1)$ in $L^2(R^d)$ as $\tau \to 0$.  Note that
$$
\int^1_0\int_{R^d} (\lambda^\tau\,\psi_t+(\lambda^\tau\nabla\phi+\lambda^\tau
(\nabla\phi^\tau-\nabla\phi)\cdot\nabla \psi)\,dxdt=
\int_{R^d} (\rho_0\psi(0)-\lambda^\tau (1)\psi(1))\,dx.
\leqno (4.17) $$
Since $\lambda^\tau$ is uniformly bounded in $L^\infty (0,1;L^1(R^d)
\cap L^\infty(R^d))$, it follows from (4.9) and (4.17) that
$\lambda^*$ is the weak solution to (4.13). Since (4.13) has the unique weak
solution, we conclude  
$\lambda^\tau$ converges weakly to $\lambda$ in $L^2((0,1)\times R^d)$ 
and $\lambda^\tau(1)$ weakly to $\lambda(1)$ in $L^2(R^d)$ as $\tau\to 0$. 
\vspace{2mm}

\noindent\underline{Step 4} Note that
$$
\int^1_0 \int_{R^d} \lambda^\tau (\psi_t+
\nabla\phi^\tau\cdot\nabla \psi)\,dxdt=
\int_{R^d}(\psi(0)\rho_0-\psi(1)\lambda^\tau(1))\,dx.
\leqno (4.18) $$
for all $\psi\in W^{1,\infty}((0,1) \times R^d)$.  Since 
$$
(\phi^\tau-\phi)_t+\nabla\phi^\tau\cdot\nabla(\phi^\tau-\phi)
-\frac{1}{2}|\nabla\phi^\tau-\nabla\phi|^2=0
$$
by setting $\psi=\phi^\tau-\phi$ in (4.18), we obtain
$$
\int_{R^d}((\phi^\tau(0)-\phi(0))\rho_0(x)-\tau\,\lambda^\tau(1)h(x))\,dx
=\frac{1}{2}\int^1_0\int_{R^d} \lambda^\tau\,|\nabla\phi^\tau-\nabla\phi|^2
\,dxdt
$$
Similarly, since
$$
(\phi^\tau-\phi)_t+\nabla\phi \cdot\nabla (\phi^\tau-\phi)
+\frac{1}{2}|\nabla\phi^\tau-\nabla\phi|^2=0,
$$ 
we have
$$
\int_{R^d}((\phi^\tau(0)-\phi(0))\rho_0(x)-\tau\,\lambda(1)h(x))\,dx
=-\frac{1}{2}\int^1_0\int_{R^d} \lambda\,|\nabla\phi^\tau -\nabla\phi|^2
\,dxdt.
$$
From (4.3) and (4.8) there exists a subsequence of
$\ds\frac{\phi^\tau(0)-\phi}{\tau}$ that converges weakly in $L^2(R^d)$
as $\tau\to 0$. 
Since $\lambda,\; \lambda^\tau \ge 0$ a.e. in $(0,1) \times R^d$
and $\lambda^\tau(1)$ converges weakly to $\lambda(1)$ as $\tau \to 0$,
we have
$$
\lim_{\tau \to 0} \int_{R^d}\frac{\phi^\tau(0)-\phi(0)}{\tau}\rho_0(x)\,dx=
\int_{R^d} \lambda(1)h(x)\,dx.
$$
Hence
$$
J^\prime(v)(h)=\int_{R^d}(\rho_1-\lambda(1))h(x)\,dx.
\leqno (4.19) $$

\noindent\underline{Step 5} Assume $\bar{v}$ attains the minimum of 
$J(v)$ in (2.8)
and $\bar{v}$ is Lipshitz and semi-convex.  
Then $J^\prime(\bar{v})(h)=0$ for all $h \in C^2_0(R^d)$
and thus from (4.19) $\bar{\lambda}(1)=\rho_1$ a.e., where $\bar{\lambda}$ 
is the weak solution to (4.13) with $\phi=\bar{\phi}$. Thus, (3.3) holds
with $\rho=\bar{\lambda}$.  Since 
$$
\phi_t+\nabla\phi\cdot\nabla\phi-\frac{1}{2}|\nabla\phi|^2=0,
$$
it follows from (4.12) with $\psi=\phi$ that 
$$
\int_{R^d} (\rho_0(x) \bar{u}(x)+\rho_1(x)\bar{v}(x))\,dx
=\int^1_0 \int_{R^d}\frac{1}{2}\, \bar{\lambda}(t,x)|\bar{V}(t,x)|^2\,dxdt.
$$
which shows (3.4).

\section{General Case}

In this section we prove (3.3)--(3.4) for the general case 
$\ds c(x-y)=\frac{1}{p}|x-y|^p,\;1<p<\infty$. Assume $v$ 
is Lipschitz and semi-convex. It the follows from \cite{E} that 
$$
\phi_t+\frac{1}{q}|\nabla \phi|^q=0,\quad \phi(1,x)=v(x)
\leqno (5.1) $$
has a unique viscosity solution $\phi \in W^{1,\infty}((0,1)\in R^d)$
satisfying (4.2).
\vspace{2mm}

\noindent\underline{Step 1}
For $h\in C_0^2(R^d)$ let $v^\tau=v+\tau\,h$.  
If for $x\in R^d$, $y^\tau=y^\tau(x) \in R^d$ 
attains the maximum  of $y \to v^\tau(y)-c(x-y)$, then
$$
x-y^\tau=-|\nabla v^\tau(y^\tau)|^{q-2}\nabla v^\tau(y^\tau)
$$
Since $v^\tau \in W^{1,\infty}(R^d)$, $|y^\tau(x)-x| \le \alpha$
for some $\alpha$ uniformly in $x$ and $\tau$.
Since as shown in Section 4
$$
\tau\,h(y)\le \pi^\tau(0,x)-\phi(0,x)\le -\tau h(y^\tau),
$$
$$
\int_{R^d}|\phi^\tau(0,x)-\phi(0,x)|\,dx
=\tau \int_{R^d} (|h(y(x))|+|h(y^\tau(x))|)\,dx.
$$
Since $h$ is compactly supported, it follws that 
there exists a constant $M$ (depends on $h$) such that
$$
\int_{R^d} |\phi^\tau(0,x)-\phi(0,x)|\,dx \le M\,\tau.
\leqno (5.2) $$
Since from (5.1)
$$
\int^1_0\int_{R^d}\frac{1}{q}\,|\nabla \phi^\tau|^q\,dxdt
=\int_{R^d} (\phi^\tau(0,x)-\phi^\tau(1,x))\,dx,
$$
we have
$$
\int^1_0\int_{R^d}(|\nabla \phi^\tau|^q-|\nabla\phi|^q)\,dx\to 0
$$
as $\tau \to 0$. Hence
$$
\int^1_0\int_{R^d}||\nabla\phi^\tau|^{\frac{q-2}{2}}\nabla\phi^\tau|^2
\,dxdt \to 
\int^1_0\int_{R^d}||\nabla\phi|^{\frac{q-2}{2}}\nabla\phi|^2\,dxdt
$$
as $\tau \to 0$ and therefore
$$
|\nabla\phi^\tau|^{\frac{q-2}{2}}\nabla\phi^\tau
\to |\nabla\phi|^{\frac{q-2}{2}}\nabla\phi
\quad\mbox{in  }L^2((0,1)\times R^d).
$$
Moreover, there exists a subsequence (denoted by the same) of $\tau$
such that $\nabla\phi^\tau(x)\to\nabla\phi(x)$ a.e. in $R^d$.
Since by the Lebesgue dominated convergence theorem
$$
|\nabla\phi^\tau|^{q-2}\nabla\phi^\tau
\to |\nabla\phi|^{q-2}\nabla\phi
\quad\mbox{in  }L^2((0,1)\times R^d).
\leqno (5.3) $$

\noindent\underline{Step 2} We assume that for $p>2$ (i.e., $q<2$)
$|\nabla v(\cdot)|^2 \ge c>0$ a.e. in $R^d$. Then it follows from (5.1) that
$|\nabla\phi(t,\cdot)|^2\ge c$ a.e. in $R^d$ for $t\in[0,1]$.
Let
$$
J_\epsilon =|\nabla (\eta_\epsilon *\phi)|^{q-2}\nabla (\eta_\epsilon *\phi) 
$$
Then,
$$
\nabla \cdot J_\epsilon= |\nabla (\eta_\epsilon *\phi)|^{q-4}
(|\nabla (\eta_\epsilon *\phi)|^2\,
\Delta (\eta_\epsilon *\phi)+
(q-2)\,\nabla(\eta_\epsilon *\phi)^t [D^2(\eta_\epsilon *\phi)] 
\nabla(\eta_\epsilon *\phi))
$$
Since from (4.2) $D^2(\eta_\epsilon *\phi)\ge -C$,
there exists a positive constant $C_q$ such that
$\nabla \cdot J_\epsilon \ge -C_q$ and thus
$$
(J_\epsilon,\nabla\psi) \le C_q\int_{R^d}\psi\,dx
$$
for $\psi \in W^{1,1}(R^d)$ and $\psi \ge 0$ a.e. in $R^d$.
It thus follows from (4.5)--(4.6) and the Lebesgue dominated 
convergence theorem that 
$$
(|\nabla\phi|^{q-2}\nabla\phi,\nabla\psi)\le C_q\int_{R^d}\psi\,dx.
\leqno (5.4) $$
for $\psi \in W^{1,1}(R^d)$ and $\psi \ge 0$ a.e. in $R^d$,
by letting $\epsilon\to0^+$.
\vspace{2mm}

\noindent\underline{Step 3}
Using the same arguments as in Step 3 in Section 4,
$$
\lambda_t+(|\nabla\phi^\tau|^{q-2}\nabla\phi^\tau\,\lambda)_x=0
$$
has the unique weak solution $\lambda^\tau \in L^\infty(0,1;L^1(R^d) 
\cap L^\infty(R^d))$, i.e.,
$$
\int^1_0 \int_{R^d} \lambda^\tau (\psi_t+
|\nabla\phi^\tau|^{q-2}\nabla\phi^\tau\cdot\nabla \psi)\,dxdt=
\int_{R^d}(\psi(0)\rho_0-\psi(1)\lambda^\tau(1))\,dx.
\leqno (5.5) $$
for all $\psi\in W^{1,\infty}((0,1) \times R^d)$.
Moreover $\lambda^\tau$ converges weakly to $\lambda$
in $L^2((0,1)\times R^d)$ and $\lambda^\tau(1)$
converges weakly to $\lambda(1)$ in $L^2(R^d)$ as $\tau\to0$. 
Since 
$$
(\phi^\tau-\phi)_t+|\nabla\phi^\tau|^{q-2}\nabla\phi^\tau\cdot
\nabla(\phi^\tau-\phi)+I_1=0
$$
where
$$ 
I_1=\frac{1}{q}|\nabla\phi^\tau|^q-\frac{1}{q}|\nabla\phi|^q
-|\nabla\phi^\tau|^{q-2}\nabla\phi^\tau\cdot(\nabla \phi^\tau-\nabla\phi)
\le 0
$$
By setting $\psi=\phi^\tau-\phi$ in (4.5), we obtain
$$
\int_{R^d}((\phi^\tau(0)-\phi(0))\rho_0(x)-\tau\,\lambda^\tau(1)h(x))\,dx
=-\int_{R^d} \lambda^\tau\,I_1\,dx 
$$
Similarly, since
$$
(\phi^\tau-\phi)_t+|\nabla\phi|^{q-2}\nabla\phi \cdot
\nabla(\phi^\tau-\phi)+I_2=0
$$
where 
$$
I_2=\frac{1}{q}|\nabla\phi^\tau|^q-\frac{1}{q}|\nabla\phi|^q
-|\nabla\phi|^{q-2}\nabla\phi\cdot(\nabla \phi^\tau-\nabla\phi)\ge 0,
$$
we have
$$
\int_{R^d}((\phi^\tau(0)-\phi(0))\rho_0(x)-\tau\,q(1)h(x))\,dx
=-\int \lambda\,I_2\,dx.
$$
Since $\lambda,\; \lambda^\tau \ge 0$ a.e. in $(0,1) \times R^d$
and $\lambda^\tau(1)$ converges weakly to $\lambda(1)$ as $\tau \to 0$,
it follws that
$$
\lim_{\tau \to 0} \int_{R^d}\frac{\phi^\tau(0)-\phi(0)}{\tau}\rho_0(x)\,dx=
\int_{R^d} \lambda(1)h(x)\,dx.
$$
Thus
$$
J^\prime(v)(h)=\int_{R^d}(\rho_1-\lambda(1))h(x)\,dx.
\leqno (5.6) $$
Assume $\bar{v}$ attains the minimum of 
$J(v)$ in (2.8)
and $\bar{v}$ is Lipshitz and semi-convex.
Then, 
$J^\prime(\bar{v})(h)=0$ for all $h\in C^2_0(R^d)$ and therrrfore
from (5.6) $\bar{\lambda}(1)=\rho_1$ where $\bar{\lambda}$ 
is the weak solution to 
$$
\lambda_t+(\lambda\,|\nabla\bar{\phi}|^{q-2}\nabla\phi)_x=0,
\quad \lambda(0)=\rho_0
\leqno (5.7) $$
Thus, (3.3) holds with $\rho=\bar{\lambda}$. Since 
$$  
\phi_t+|\nabla\phi|^{q-2}\nabla\phi\cdot\nabla\phi-
\frac{1}{p}|\nabla\phi|^q=0,
$$
it follows from (5.5) with $\psi=\phi$ that 
$$
\int_{R^d} (\rho_0(x) \bar{u}(x)+\rho_1(x)\bar{v}(x))\,dx
=\int^1_0 \int_{R^d}\frac{1}{p}\,\bar{\lambda}(t,x)|\bar{V}(t,x)|^p\,dxdt.
$$
which shows (3.4).

\end{document}